# Regional Boundary Gradient Strategic Sensors Characterizations

Raheam A. Al-Saphory[1*] Naseif J. Al-Jawari[2] Asmaa N. Al-Janabi[3]

1.Department of Mathematics, College of Education for Pure Sciences, Tikrit University, Tikrit, IRAQ.

2.Department of Mathematics, College of Science, Al-Mustansiriyah University, Baghdad, IRAQ.

3.Department of Mathematics, College of Science, Al-Mustansiriyah University, Baghdad, IRAQ.

* E-mail of the corresponding author: saphory@hotmail.com

**Abstract**

The purpose of this paper is to characterize regional boundary gradient strategic sensors notions for different cases of regional boundary gradient observation to be achieved. Then, the characterizations based on how to a cross from internal gradient region to the boundary gradient region. Thus, the obtained results are applied in two dimensional linear infinite distributed systems in Hilbert space where the dynamics are governed by strongly continuous semi-group. Moreover, we give the relation between the regional gradient observability on a subregion $\omega$ of the spatial domain $\Omega$ with the regional boundary gradient observability on a subregion $\Gamma$ of the boundary $\partial\Omega$ of $\Omega$. Finally, sufficient conditions of regional boundary gradient strategic sensors notions are explored, analyzed and discussed in connection with the regional boundary gradient of exact (weak) observability, positive definite observability operator and rank conditions.

## 1. Introduction

One of the first steps in engineering designing is to identify its physical part. The physical part can be belong to a large array of system classes (Curtain & Pritchard 1978). So, we focus our interest on distributed parameter systems whose dynamics can be involves partial differential equations where the states depend not on time only but also on spatial variables (Curtain & Zwart 1995). The analysis of distributed parameter systems refers to a many concepts such as controllability and observability (El Jai & Amouroux 1988). The study of these concepts are related to the notions of sensors and actuators where the characterization of sensors and actuators are playing a fundamental role in the understanding of any real systems because they are intermediates between a system and it's environment (El Jai & Pritchard 1988). The regional analysis is one of the most important notion of system theory, is focused on a state observation on a sub-region ω of the spatial domain $\Omega$. These concepts are introduced and developed by El Jai et. al. as in ref.s (El Jai *et al.* 1994), (El Jai & El Yacoubi 1993), (Zerrik 1993), (Bourray *et al.* 2014) (Bourray *et al.* 2014), (Al-Saphory *et al.* ), (R. Al-Saphory & El Jai 2002) & (Al-Saphory & .El Jai 2001) and it has been extended to the case where the sub-region $\Gamma$ is a part of the boundary $\partial\Omega$ of $\Omega$ in (Zerrik & Badraoui 2000) & (Zerrik *et al.* 1999). In the same direction one may be concerned with the regional gradient observability for a diffusion system where one is interested in knowledge of the state gradient only in a critical sub-region of the system domain without the knowledge of the state itself, this concepts has been introduced by (Zerrik & Bourray 2003). Thus, the concept of regional strategic sensors which was introduced in (Al-Saphory & Al-Joubory 2013) gives a characterization of a regional strategic sensor to achieve of regional observability. In addition, the result in (Al-Saphory *et al*. 2015) has been extended to the regional gradient strategic sensors for various types of regional gradient observability. The introduction of this concept is motivated by real situations. This is the case, for example, of energy exchange problem, where the aim is to determine the energy exchanged between a casting plasma on a plane target which is perpendicular to the direction of the flow from measurements carried out by internal thermocouples (Figure 1).

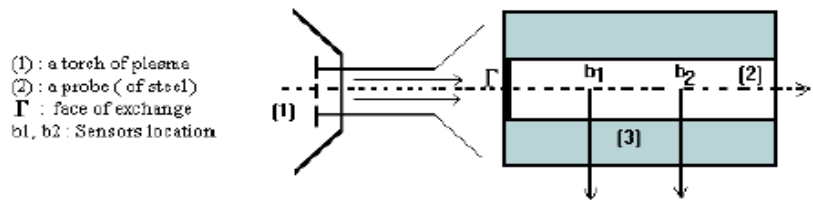

(1) : a torch of plasma
(2) : a probe ( of steel)
$\Gamma$ : face of exchange
b1, b2 : Sensors location

Figure1. The Estimation Profile of Energy Exchanged on $\Gamma$.

More precisely let the considered system with suitable state space, and suppose that the initial state $y_0$ and its gradient $\nabla y_0$ are unknowns. Suppose now that measurements are given by means of an output function (depending on the number and the structure of the sensors) (El Jai & Pritchard 1988). The problem concerns the reconstruction of the initial gradient $\nabla y_0$ on the subregion $\Gamma$ of the system domain boundary $\partial\Omega$.

The main objective of this paper is to extended these results in (Al-Saphory *et al*. 2015) to the case of regional boundary gradient strategic sensors for different cases of regional boundary gradient observation. More precisely, we discuss and analyze the relation between the regional gradient strategic sensors and the regional boundary gradient strategic sensors. This paper which organized as follows:

Second section is present the problem statement and basic definitions with characterizations on the regional boundary gradient observability. Third section is devoted to the mathematical concepts of regional boundary gradient strategic sensors in various situations are studied and we develop an approach to cross internal region to the boundary region. In the last section we gives an application about different sensors locations.

## 2. Regional Boundary Gradient Observability

In this section, we present some notions and preliminary material as in (Zerrik & Bourray 2003).

### 2.1 Problem Statement

Let $\Omega$ be a regular bounded open subset of $R^n$, with boundary $\partial\Omega$ and $\Gamma$ be sub-region of $\partial\Omega$, [0,T], $T>0$ be a time measurement interval. We denoted $Q = \Omega \times ]0,T[$, $\Sigma = \partial\Omega \times ]0,T[$. We considered distributed parabolic system is described by the following equation:

$$\begin{cases} \frac{\partial x}{\partial t}(\xi, t) = Ax(\xi, t) + Bu(t), & Q \\ x(\xi, 0) = x_0(\xi), & \Omega \\ x(\eta, t) = 0, & \Sigma \end{cases} \tag{1}$$

with the measurement is given by the output function

$$y(., t) = Cx(., t). \tag{2}$$

where $A = \sum_{i,j=1}^{n} \frac{\partial}{\partial x_j}\left(a_{ij}\frac{\partial}{\partial x_j}\right)$   with   $a_{ij} \in \mathcal{D}(\bar{Q})$  (domain of  $\bar{Q}$) is a second order linear differential operator, which is generated a strongly continuous semi-group $\left(S_A(t)\right)_{t \geq 0}$ on the Hilbert space  $X$  and is self-adjoint with compact resolvent. Suppose that  $-A$  is elliptic, i.e., there exits  $\alpha > 0$  such that $\sum_{i,j=1}^{n} a_{ij}\xi_i\xi_j \geq \alpha \sum_{j=1}^{n}\left|\xi_j\right|^2$,  almost  everywhere  (a.e.)  on  Q ,  $\forall \xi = (\xi_1, \dots, \xi_n) \in R^n$.  The  operator $B \in L(R^p, X)$  and  $C \in L(X, R^q)$,  depend on the structure of actuators and sensors (El Jai & Amouroux 1988). The space  $X, U$ and $\mathcal{O}$  be separable Hilbert spaces where  $X$  is the state space,  $U = L^2(0, T, R^p)$  is the control space and  $\mathcal{O} = L^2(0, T, R^q)$  is the observation space where  $p$  and  $q$  are the numbers of actuators and sensors. Under the given assumption, the system (1) has a unique solution (Curtain & Pritchard 1978) & (Curtain & Zwart 1995) given as follows:

$$x(\xi, t) = S_A(t)x_0(\xi) + \int_0^t S_A(t-\tau)Bu(\tau)d\tau \tag{3}$$

The problem is that, how to give sufficient conditions for regional boundary gradient strategic sensors. These conditions which enable to achieve the regional boundary gradient observability with different type of strategic sensors depending on the domain region, dimensional systems and structure sensors.

### 2.2 Definition and Characterizations

In this sub-section we are recall and introduced some basic concepts about the regional boundary gradient observability and sensors structures. For this purpose, we give the linear autonomous system of (1)

$$\begin{cases} \frac{\partial x}{\partial t}(\xi, t) = Ax(\xi, t), & Q \\ x(\xi, 0) = x_0(\xi), & \Omega \\ x(\eta, t) = 0. & \Sigma \end{cases} \tag{4}$$

The solution of (4) is given by the following form

$$x(\xi, t) = S_A(t)x_0(\xi), \quad \forall t \in [0, T] \tag{5}$$

The measurements are obtained through the output function by using of zone, pointwise which may located in  $\Omega$  (or  $\partial\Omega$)  given by the following form ((El Jai & Pritchard 1988):

$$y(., t) = Cx(\xi, t) \tag{6}$$

• We first recall a sensor is defined by any couple  $(D, f)$, where  $D$  is spatial support represented by a nonempty part of  $\bar{\Omega}$  and f  represents the distribution of the sensing measurements on  $D$.

Depending on the nature of $D$ and $f$, we could have various type of sensors. A sensor may be pointwise if $D = \{b\}$ with $b \in \overline{\Omega}$ and $f = \delta(.-b)$, where $\delta$ is the Dirac mass concentrated at $b$. In this case the operator C is unbounded and the output function (2) can be written in the following form

$$y(t) = x(b, t)$$

- A sensor may be zonal when $D \subset \overline{\Omega}$ and $f \in L^2(D)$. In this case the output function (2) can be written by the form

$$y(t) = \int_D x(\xi, t) f(\xi) d\xi$$

- In the case of boundary zone sensor, we consider $D_i = \Gamma_i \subset \partial\Omega$ and $f_i \in L^2(\Gamma_i)$, the output function (2) can be written as

$$y(., t) = Cx(., t) = \int_{\Gamma_i} x(\eta, t) f_i(\eta) d\eta$$

- Now, we define the operator

$$K: x \in X \to Kx = CS_A(.)x \epsilon \, \mathcal{O}$$

we note that $K^*: \mathcal{O} \to X$ is the adjoint operator of $K$ defined by

$$K^* y^* = \int_0^t S_A^*(s) C^* y^*(s) ds$$

- Consider the operator

$$\nabla: \begin{cases} H^1(\Omega) \to (H^1(\Omega))^n \\ x \to \nabla x = (\frac{\partial x}{\partial \xi_1}, \dots, \frac{\partial x}{\partial \xi_n}) \end{cases}$$

and the adjoint of $\nabla$ denotes by $\nabla^*$ is given by

$$\nabla^*: \begin{cases} (H^1(\Omega))^n \to H^1(\Omega) \\ x \to \nabla^* x = v \end{cases}$$

where $v$ is a solution of the Dirichlet problem

$$\begin{cases} \Delta v = -div(x) & in \; \Omega \\ v = 0 & in \; \partial\Omega \end{cases}$$

- The trace operator of order zero is given by

$$\gamma_0: H^1(\Omega) \to H^{1/2}(\partial\Omega)$$

Thus, the extension of the trace operator of order zero which is denoted by $\gamma$ is linear, subjective and continuous, defined as an operator $\gamma: (H^1(\Omega))^n \to (H^{1/2}(\partial\Omega))^n$ and the adjoints are respectively given by $\gamma_0^*, \gamma^*$.

- For $\Gamma \subset \partial\Omega$, we consider a gradient restriction operator $\chi_\Gamma : (H^{1/2}(\partial\Omega))^n \to H^{1/2}(\Gamma))^n$, and

$\tilde{\chi}_\Gamma : H^{\frac{1}{2}}(\partial\Omega) \to H^{\frac{1}{2}}(\Gamma)$,

where the adjoint are respectively given by $\chi_\Gamma^*$, $\tilde{\chi}_\Gamma^*$.

- Also, if $\omega \subset \Omega$ we consider the operator

$$\chi_\omega : \begin{cases} (H^1(\Omega))^n \to (H^1(\omega))^n \\ x \quad \to \quad \chi_\omega x = x \mid_\omega \end{cases}$$

It's adjoint is denoted by $\chi_\omega^*$

- Finally, we introduced the operator $H = \chi_\Gamma \gamma \nabla K^*$ from $\mathcal{O}$ into $(H^{1/2}(\Gamma))^n$ Now, let us denoted the system (4) together with the output (6) by (4)-(6).

**Definition 2.1:** The systems (4)-(6) are said to be exactly regionally gradient observable on $\omega$ (exactly $\omega_G$-observable), if

$$Im\chi_\omega \nabla K^* = (H^1(\omega))^n$$

**Definition 2.2:** The systems (4)-(6) are said to be weakly regionally gradient observable on $\omega$ (weakly $\omega_G$-observable), if

$$\overline{Im\chi_\omega \nabla K^*} = (H^1(\omega))^n$$

**Definition 2.3:** The systems (4)-(6) are said to be exactly regionally boundary gradient observable on $\Gamma$ (exactly $\Gamma_G$ –observable) if

$$Im\, H = Im\chi_\Gamma \gamma \nabla K^* = (H^{1/2}(\Gamma))^n$$

**Definition 2.4:** The systems (4)–(6) are said to be weakly regionally boundary gradient observable on $\Gamma$ (weakly $\Gamma_G$ –observable) if

$$\overline{ImH} = \overline{Im\chi_\Gamma \gamma \nabla K^*} = (H^{1/2}(\Gamma))^n$$

**Remark 2.5:** The definition (2.4) is equivalent to say that the systems (4)-(6) are weakly $\Gamma_G$ –observable if

$$ker\, H^* = ker\, K\nabla^* \gamma^* \chi_\Gamma^* = \{0\}$$

**Proposition 2.6:** The systems (4) -(6) are exactly $\Gamma_G$-observable if and only if there exists $c > 0$ such that, for all $x^* \in (H^{1/2}(\Gamma))^n$,

$$\| \chi_\Gamma x^* \|_{(H^{1/2}(\Gamma))^n} \leq c \| K\nabla^* \gamma^* \chi_\Gamma^* x^* \|_{\mathcal{O}} \qquad (7)$$

***Proof:*** The proof of this property is deduced from the following general result (Curtain & Pritchard 1978).

Let $E, F, G$ be reflexive Banach space and $f \in L(E, G)$, $g \in L(F, G)$, then the following properties are equivalent:

a. $Im\, f \subset Im\, g$.

b. There exists $c > 0$ such that $\| f^* x^* \|_{E^*} \leq c \| g^* x^* \|_{F^*}$, $\forall x^* \in G^*$.

If we apply this result, considered $E = G = (H^{1/2}(\Gamma))^n, F = \mathcal{O}, f = Id_{(H^{1/2}(\Gamma))^n}$ and $g = \chi_\Gamma \gamma \nabla K^*$. We obtain the considered inequality (7).□

**Remark 2.7:** From the previous proposition 2.6, we can get the following results:

1. The regional state reconstruction will be more precise than the boundary of the domain $\partial\Omega$ if we estimate the state in whole the boundary domain.

2. From (7) there exists a reconstruction error operator that gives the estimation $\tilde{x}_0$ of the initial state $x_0$ in $\Gamma$,   and if we put $e = x_0 - \tilde{x}_0$ , we have

$$\|e\|_{H^{1/2}(\Gamma)} \leq \|e\|_{H^{1/2}(\partial\Omega)}$$

and then

$$\|x_0 - \tilde{x}_0\|_{H^{1/2}(\Gamma)} \leq \|x_0 - \tilde{x}_0\|_{H^{1/2}(\partial\Omega)}$$

where, $x_0$ is the exact state of the system and $\tilde{x}_0$ is the estimated state of the system.

3. If a system is exactly observable in $\partial\Omega$ then it is regionally exactly observable in $\Gamma$.

Now, we can deduced that :

**Proposition 2.8:** If a system is exactly regionally boundary observable (exactly $\Gamma_B$-observable) then it is exactly $\Gamma_G$-observable.

***Proof:*** Since the system is exactly $\Gamma_B$-observable then there exists $\gamma > 0$ such that $\forall x_0 \in L^2(\Gamma)$, we have,

$$\|x_0\|_{L^2(\Gamma)} \leq \gamma \|K\gamma_0^* \tilde{\chi}_\Gamma^* x_0\|_{L^2(0,T,\mathcal{O})}, \ \forall \gamma > 0$$

since $(L^2(\Gamma))^n \subset L^2(\Gamma)$, then

$$\|\nabla x_0\|_{(L^2(\partial\Omega))^n} = \|x_0\|_{(L^2(\Gamma))^n} \leq \|x_0\|_{L^2(\Gamma)} \ , \ \forall x_0 \in L^2(\Gamma) \text{ where,}$$

$$L^2(\Gamma) = \{x_0 : \int_\Gamma |x_0|^2 < \infty\} \text{ and}$$

$$(L^2(\Gamma))^n = \{\nabla x_0 = g_i : \int_\Gamma |g_i|^2 < \infty, \ g_i = \frac{\partial x_0}{\partial \xi_i} \ \forall i = 1,2,\dots\}. \tag{8}$$

To prove that $\|x_0\|_{(L^2(\Gamma))^n} \leq c \|K\nabla^* \gamma^* \chi_\Gamma^* x_0\|_{L^2(0,T,\mathcal{O})}$. Thus, from (8) and since a system is exactly $\Gamma_B$-observable, then

there exists $\gamma > 0$ and $c > 0$ which allow to choose $\gamma = \frac{1}{c}$ where $c$ is given by

$$c = \frac{\|K\gamma_0^* \tilde{\chi}_\Gamma^* x_0\|_{\mathcal{O}}}{\|K\nabla^* \gamma^* \chi_\Gamma^* x_0\|_{\mathcal{O}}} \tag{9}$$

and then

$$\|x_0\|_{(L^2(\Gamma))^n} \leq \|x_0\|_{L^2(\Gamma)} \leq \gamma \|K\gamma_0^* \tilde{\chi}_\Gamma^* x_0\|_{\mathcal{O}} \tag{10}$$

substitute (9) in (10), we obtain

$$\|x_0\|_{(L^2(\Gamma))^n} \leq \|K\nabla^* \gamma^* \chi_\Gamma^* x_0\|_{\mathcal{O}}.$$

Therefor, this system is exactly $\Gamma_G$-observable   with $\gamma = 1$.□

## 3. Boundary Gradient Strategic Sensors

The purpose of this section is to link the regional boundary gradient observability notion with the sensors structure. Consider now the system (4) observed by $q$ sensors $(D_i, f_i)_{1 \leq i \leq q}$, which may be pointwise or zonal.

### 3.1 $\Gamma_G$-Strategic Sensors

**Definition 3.1:** A sensor $(D, f)$ is boundary gradient strategic on $\Gamma$ ($\Gamma_G$-strategic) if the observed system is weakly $\Gamma_G$-observable.

**Definition 3.2:** A suite of sensors $(D_i, f_i)_{1 \leq i \leq q}$ are boundary gradient strategic on $\Gamma$ if there exists at least one sensor $(D_1, f_1)$ which is weakly $\Gamma_G$-observable.

**Corollary 3.3:** A sensor is $\Gamma_G$-strategic if the corresponding observed system is exactly $\Gamma_G$ –observable.

***Proof:*** Since the system is exactly $\Gamma_G$ –observable, then we have

$$Im\, H = Im\chi_\Gamma \gamma \nabla K^* = (H^{1/2}(\Gamma))^n.$$

From decomposition subspace of direct sum in Hilbert space, we represent $(H^{1/2}(\partial\Omega))^n$ by the unique form [6]

$$ker\chi_\Gamma + Im\chi_\Gamma^* \chi_\Gamma \gamma \nabla K^* = (H^{1/2}(\partial\Omega))^n$$

we obtain

$$ker\, K(t)\, \nabla^* \gamma^* \chi_\Gamma^* = \{0\}$$

this is equivalent to

$$\overline{Im\, \chi_\Gamma \gamma \nabla K^*} = (H^{1/2}(\Gamma))^n.$$

Finally, we can deduce this system is weakly $\Gamma_G$–observable and therefore this sensor is $\Gamma_G$-strategic.□

**Corollary 3.4:** A sensor is $\Gamma_G$-strategic if and only if the operator $N_\Gamma = HH^*$ is positive definite.

***Proof:*** Since a sensor is $\Gamma_G$-strategic this means that a system is weakly $\Gamma_G$- observable. Let $x^* \in (H^{\frac{1}{2}}(\Gamma))^n$ such that

$$< N_\omega x^*, x^* >_{(H^{1/2}(\Gamma))^n} = 0 \quad \text{then} \quad \|H^* x^*\|_{\mathcal{O}} = 0$$

and therefore $x^* \in ker H^*$, thus, $x^* = 0$, i.e., $N_\omega$ is positive definite.

Conversely, let $x^* \in (H^{1/2}(\Gamma))^n$ such that

$$H^* x^* = 0, \text{ then } < H^* x^*, H^* x^* >_{\mathcal{O}} = 0$$

and thus,

$$< N_\omega x^*, x^* >_{(H^{1/2}(\Gamma))^n} = 0.$$

Hence $x^* = 0$ thus the system is weakly $\Gamma_G$- observable and therefore a sensor is $\Gamma_G$-strategic.□

Thus, from previous corollaries we can obtain the following results:

**Remark 3.5:** We can deduce that:

1. If the system is exactly $\Gamma_G$–observable then the system is weakly $\Gamma_G$–observable and therefore the sensor is $\Gamma_G$-strategic.

2. A sensor which is regional boundary gradient strategic in $\Gamma_G^1$ ($\Gamma_G^1$-strategic sensor) for a system where $\Gamma_G^1 \subset \Gamma_G$, is regional boundary gradient strategic sensor in $\Gamma_G^2$ ($\Gamma_G^2$-strategic sensor) for any $\Gamma_G^2 \subset \Gamma_G^1$.

### 3.2 Crossing Approach from Internal to Boundary Case

In this approach we deal with the regional boundary gradient strategic sensors in $\Gamma$ ($\Gamma_G$-strategic sensor) as an internal regional gradient strategic sensors. In this case, we introduced the following extension operator (Zerrik & Bourray 2003)

- Let $\Re : (H^{1/2}(\partial\Omega))^n \to (H^1(\Omega))^n$, which is continuous and linear defined by

$$\gamma \Re h(\xi,t) = h(\xi,t), \forall h(\xi,t) \in (H^{1/2}(\partial\Omega))^n \tag{11}$$

- Let $E = \bigcup_{x \in \Gamma} B(x,r)$ and $\overline{\omega}_r = E \cap \Omega$, where $B(x,r)$ is a ball of radius $r$, ($r > 0$ is an arbitrary and sufficiently small real) and centered in $x(\xi,r)$ and $\Gamma$ is a part of $\overline{\omega}_r$ (Figure 2).

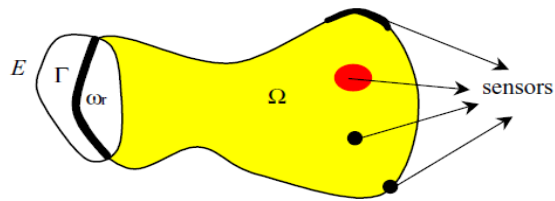

Figure 2. The Domain $\Omega$, Subregion $\overline{\omega}_r$ and the Region $\Gamma$.

### 3.3 $\Gamma_G$-strategic sensors and $\overline{\omega}_{r_G}$-strategic sensors

In this subsection we explain the link between the $\Gamma_G$-strategic sensor and the regionally gradient strategic sensor on $\overline{\omega}_r$ ($\overline{\omega}_{r_G}$-strategic sensor).

**Remark 3.6:** From the above results we link the internal regional gradient observability in $\overline{\omega}_r$ to the boundary gradient case on $\Gamma$, so we can deduced the following corollary.

**Corollary 3.7:** If a sensor is $\overline{\omega}_{r_G}$– strategic then it is $\Gamma_G$-strategic.

**Proof:** Since the sensor is $\overline{\omega}_{r_G}$– strategic in $\overline{\omega}_r$, this mean that the system is weakly $\overline{\omega}_{r_G}$– observable in $\overline{\omega}_r$.

Thus, the system is weakly $\Gamma_G$– observable (Zerrik & Bourray 2003)

Therefor, the sensor is $\Gamma_G$-strategic.□

**Corollary 3.8:** A sensor is $\Gamma_G$-strategic if the system is exactly $\overline{\omega}_{r_G}$–observable.

**Proof:** Let $x(\xi,t) \in (H^{1/2}(\Gamma))^n$ and $\bar{x}(\xi,t)$ be an extension to $(H^{1/2}(\partial\Omega))^n$. By using equation (11) and trace theorem there exists $\Re\bar{x}(\xi,t) \in (H^1(\Omega))^n$, with bounded support such that

$$\gamma \Re\bar{x}(\xi,t) = \bar{x}(\xi,t).$$

Since the system is exactly $\overline{\omega}_{r_G}$–observable, then the system is weakly $\overline{\omega}_{r_G}$–observable (Al-Saphory *et al.* 2015).

since a system is weakly $\overline{\omega}_{r_G}$–observable then a system is weakly $\Gamma_G$–observable (Al-Saphory R. (2002) & ( Zerrik & Bourray 2003) . Thus, a sensor is $\Gamma_G$-strategic.□

**Remark 3.9:** We can deduce that:

1. If the system is exactly $\overline{\omega}_{r_G}$–observable then it is exactly $\Gamma_G$–observable.

i.e., there exists an operator $\chi_{\bar{\omega}_r} \nabla K^*: \mathcal{O} \to (H^1(\omega_r))^n$ given by

$$H_{\bar{\omega}_r} y(.,t) = \chi_{\bar{\omega}_r} \nabla K^* y(.,t) = \chi_{\bar{\omega}_r} \Re \tilde{x}(\xi, t).$$

Hence,

$$\chi_\Gamma \left( \gamma \chi_{\bar{\omega}_r} \nabla K^* y(.,t) \right) = x(\xi, t).$$

Where $x(\xi,t) \in (H^{1/2}(\Gamma))^n$ and $\bar{x}(\xi,t)$ be an extension to $(H^{1/2}(\partial\Omega))^n$.

2. If the system is weakly $\bar{\omega}_{r_G}$–observable then it is weakly $\Gamma_G$–observable.

3. An extension of these results can be applied for different cases of regional exponential general observability (Al-Saphory & Jaafar 2015). And, to the regional exponential reduced observability in distributed parameter systems (Al-Saphory & Al-Mullah 2015 ).

The concept of boundary gradient strategic on $\Gamma$ can be characterized by the following main result may be called rank conditions:

assume that there exists a complete set of eigenfunctions $(\varphi_{nj})_{n \in I \ \& \ j=1,\ldots,r_n}$ of $A$ in $H^1(\Omega)$ associated with eigenvalue $\lambda_n$ of multiplicities $r_n$ and $r_n = sup_{n \in I} r_n$ is finite. For $\bar{x} = (x_1,\ldots,x_{n-1})$ and $\bar{n} = (n_1,\ldots,n_{n-1})$. Suppose that the function $\psi_{\bar{n}j}(\bar{x}) = \chi_\Gamma \gamma \nabla \varphi_{nj}(x)$, $n \in I$, form a complete set in $(H^{1/2}(\Gamma))^n$.

**Theorem 3.10:** Assume that $sup \ r_n = r < \infty$, then the suite of sensors $(D_i, f_i)_{1 \le i \le q}$, $\Gamma_G$-strategic sensor if and only if

$q \ge r$, $rank \ G_n = r_n \ \ \forall n \ge 1$, $where \ G_n = (G_n)_{ij} \ for \ 1 \le i \le q, 1 \le j \le r_n$

and

$$(G_n)_{ij} = \begin{cases} \sum_{k=1}^m \frac{\partial \psi_{\bar{n}j}}{\partial x_k}(b_i) & \text{in the pointwise case} \\ \sum_{k=1}^n <\frac{\partial \psi_{\bar{n}j}}{\partial x_k}, f_i>_{L^2(D_i)} & \text{in the zonal case} \end{cases}$$

**Proof:** First, we know that if a system is weakly $\Gamma_G$ –observable then is equivalent to $[K\nabla^* \gamma^* X_\Gamma^* x^* = 0 \Longrightarrow x^* = 0]$ which allows to say that the sequence of sensors $(D_i, f_i)_{1 \le i \le q}$ is $\Gamma_G$-strategic if and only if

$$\left\{ x^* \in (H^{1/2}(\Gamma))^n \middle| <Hy, x^*>_{(H^{1/2}(\Gamma))^n} = 0, \forall y \in \mathcal{O} \right\} = \{0\}$$

Suppose that the suite of sensors $(D_i, f_i)_{1 \le i \le q}$ is $\Gamma_G$-strategic on $\Gamma$, but for a certain $n \in N$, $rank \ G_n \ne r_n$, then there exists a vector $x_n = (x_{n_1}, x_{n_2}, \ldots, x_{n_{r_n}})^{tr} \ne 0$, such that $G_n x_n = 0$.

So, we can construct a nonzero $x_0 \in H^{1/2}(\Gamma)$ considering $<x_0, \psi_{pj}>_{H^{1/2}(\Gamma)} = 0$ if $p \ne n$, and $<x_0, \psi_{nj}>_{H^{1/2}(\Gamma)} = x_{nj}$, $1 \le j \le r_n$.

Let $x_0 = \sum_{j=1}^{r_n} x_{nj} \psi_{nj}$, $X_0 = (x_0, x_0, \ldots, x_0)$ then:

$$<Hy, x_0>_{(H^{1/2}(\Gamma))^n} = \sum_{k=1}^n <\tilde{\chi}_\Gamma \gamma_0 \frac{\partial}{\partial \xi_k}(K^* y), \tilde{\chi}_\Gamma^* x_0>_{H^{1/2}(\Gamma)}$$

$$= \sum_{k=1}^n <\frac{\partial}{\partial \xi_k}(\tilde{x}(T)), \gamma_0 \tilde{\chi}_\Gamma^* x_0>_{L^2(\partial\Omega)}$$

where $\tilde{x}$ is the solution of the following system:

$$\begin{cases} \frac{\partial \tilde{x}}{\partial t}(\xi,t) = A^* \tilde{x}(\xi,t) + \sum_{i=1}^{q} f_i y_i (T-t) & Q \\ \tilde{x}(\xi,0) = 0 & \Omega \\ \tilde{x}(\eta,t) = 0 & \Sigma \end{cases} \qquad (12)$$

Consider the system:

$$\begin{cases} \frac{\partial \psi}{\partial t}(\xi,t) = -A\varphi(\xi,t) & Q \\ \psi(\xi,0) = \gamma_0^* \tilde{\chi}_{\Gamma}^* x_0 & \Omega \\ \psi(\eta,t) = 0 & \Sigma \end{cases} \qquad (13)$$

multiply (12) by $\frac{\partial \psi}{\partial \xi_k}$ and integrate on $Q$, we obtain:

$$\int_Q \frac{\partial \psi}{\partial \xi_k}(\xi,t) \frac{\partial \tilde{x}}{\partial t}(\xi,t) d\xi \, dt = \int_Q A^* \tilde{x}(\xi,t) \frac{\partial \psi}{\partial \xi_k}(\xi,t) d\xi \, dt$$

$$+ \int_Q \left( \sum_{i=1}^{q} \delta_{b_i} y_i (T-t) \right) \frac{\partial \psi}{\partial \xi_k}(\xi,t) d\xi dt$$

but we have

$$\int_Q \frac{\partial \varphi}{\partial \xi_k}(\xi,t) \frac{\partial \tilde{x}}{\partial t}(\xi,t) d\xi \, dt = \int_{\partial \Omega} \left[ \frac{\partial \varphi}{\partial \xi_k}(\xi,t) \tilde{x}(\xi,t) d\xi \right]_0^T + \int_Q A \frac{\partial \psi}{\partial \xi_k}(\xi,t) \tilde{x}(\xi,t) d\xi \, dt$$

$$= \int_{\partial \Omega} \frac{\partial \psi}{\partial \xi_k}(\xi,t) \tilde{x}(\xi,t) d\xi \; + \; \int_Q A \frac{\partial \psi}{\partial \xi_k}(\xi,t) \tilde{x}(\xi,t) d\xi \, dt$$

then

$$\int_{\partial \Omega} \frac{\partial \psi}{\partial \xi_k}(\xi,t) \tilde{x}(\xi,t) d\xi = -\int_Q A \frac{\partial \psi}{\partial \xi_k}(\xi,t) \tilde{x}(\xi,t) d\xi + \int_Q A^* \tilde{x}(\xi,t) \frac{\partial \psi}{\partial \xi_k}(\xi,t) d\xi \, dt \quad + \int_Q \left( \sum_{i=1}^{q} \delta_{b_i} y_i (T-t) \right) \frac{\partial \psi}{\xi_k}(\xi,t) d\xi dt.$$

Integrating by parts we obtain

$$\int_{\partial \Omega} \frac{\partial \psi}{\partial \xi_k}(\xi,t) \tilde{x}(\xi,t) d\xi \quad = \quad -\int_\pi \frac{\partial \tilde{x}(\eta,t)}{\partial \nu_{A^*}} \frac{\partial \psi}{\partial \xi_k}(\eta,t) d\eta dt \quad + \quad \int_\pi \frac{\partial}{\partial \nu_{A^*}} \left( \frac{\partial \psi}{\partial \xi_k}(\eta,t) d\eta dt \right) \tilde{x}(\eta,t) d\eta dt$$

$$+ \int_Q \left( \sum_{i=1}^{q} \delta_{b_i} y_i (T-t) \right) \frac{\partial \psi}{\partial \xi_k}(\xi,t) d\xi dt$$

the boundary conditions give

$$\int_{\partial \Omega} \frac{\partial \psi}{\partial \xi_k}(\xi,t) \tilde{x}(\xi,t) d\xi = \int_Q \left( \sum_{i=1}^{q} \delta_{b_i} y_i (T-t) \right) \frac{\partial \psi}{\partial \xi_k}(\xi,t) d\xi dt.$$

Thus,

$$\int_{\partial \Omega} \psi(\xi,t) \frac{\partial \tilde{x}}{\partial \xi_k}(\xi,T) d\xi = -\sum_{i=1}^{q} \int_0^T \frac{\partial \psi}{\partial \xi_k}(b_i,t) \, y_i (T-t) dt$$

and we have

$$< \chi_\Gamma \gamma \nabla K^*, x_0 >_{(H^{1/2}(\Gamma))^n} = \sum_{k=1}^{n} \int_\Omega \frac{\partial \tilde{x}}{\partial \xi_k}(\xi,t) \psi(\xi,t) d\xi$$

$$= -\sum_{K=1}^{q}\int_0^T\sum_{k=1}^n\frac{\partial\psi}{\partial\xi_k}(b_i,t)\,y_i(T-t)dt$$

but

$$\psi(\xi,t) = \sum_{p=1}^{\infty}e^{-\lambda_p(T-t)}\sum_{j=1}^{r_p}<x_0,\psi_{pj}>_{L^2(\omega)}\psi_{pj}.$$

Then,

$$\sum_{k=1}^n\frac{\partial\psi}{\partial\xi_k}(b_i,t) = \sum_{p=1}^{\infty}e^{-\lambda_p(T-t)}\sum_{j=1}^{r_p}<x_0,\varphi_{pj}>_{L^2(\omega)}\sum_{k=1}^n\frac{\partial\psi}{\partial\xi_k}(b_i)$$

$$= \sum_{p=1}^{\infty}e^{\lambda_p(T-t)}(G_px_p)_i$$

therefore,

$$<\chi_\Gamma\gamma\nabla K^*y,x_0>_{(H^{1/2}(\Gamma))^n} = -\sum_{K=1}^{q}\int_0^T\sum_{p=1}^{\infty}e^{\lambda_p(T-t)}(G_px_p)\,y_i(T-t)dt$$

thus,

$$<\chi_\Gamma\gamma\nabla K^*y,x_0>_{(H^{1/2}(\Gamma))^n} = -\sum_{i=1}^{q}\int_0^T e^{\lambda_n(T-t)}(G_nx_n)_i\,y_i(T-t)dt = 0.$$

This is true for all $\quad y \in L^2(o,T:R^q)$, then $\quad X_0 \in ker\,H^*$ which contradicts the assumption that the suite of sensor is $\Gamma_G$-strategic.□

**Corollary 3.11:** If the systems (4)-(6) are exactly $\Gamma_G$ –observable, then the rank condition in theorem (3.10) is satisfied and the sensor is $\Gamma_G$-strategic.

**Remark 3.12:** From the above   results, we can deduce that:

1. The theorem 3.10 implies that the required number of sensors is greater than or equal to the largest multiplicity of the eigenvalues.

2. By infinitesimally deforming the domain, the multiplicity can be reduced to one (El Jai & El Yacoubi 1993). Consequently,   $\Gamma_G$-observability can be achieved using only one sensor.

3. We can show that various sensors which are not gradient strategic in usual sense for the system but may be $\Gamma_G$-strategic and achieve the $\Gamma_G$-observability as in (Al-Saphory 2002)

## 4. Applications to Sensors Locations

In this section, we explore various results related to different types of measurements and we consider a two dimensional diffusion system defined on $\Omega = ]0,a_1[\times]0,a_2[$ by

$$\begin{cases}\frac{\partial x}{\partial t}(\xi_1,\xi_2,t) = \frac{\partial^2 x}{\partial\xi_1^2}(\xi_1,\xi_2,t) + \frac{\partial^2 x}{\partial\xi_2^2}(\xi_1,\xi_2,t) & Q\\ x(\xi_1,\xi_2,0) = x_0(\xi_1,\xi_2) & \Omega\\ x(\xi,\eta,t) = 0 & \Sigma\end{cases} \qquad (16)$$

Let $\Gamma = ]0,a_1[\times\{a_2\}$, the eigenfunctions and the eigenvalues of the system (16) are given by:

$$\varphi_{inm}(\xi_1,\xi_2) = \frac{2}{\sqrt{a_1a_2}}sin\,n\pi\frac{\xi_1}{a_1}sin\,m\pi\frac{\xi_2}{a_2} \qquad (17)$$

associated with eigenvalue

$$\lambda_{nm} = -(\frac{i^2}{a_1^2} + \frac{j^2}{a_2^2})\pi^2 \qquad (18)$$

if we suppose that $\frac{a_1^2}{a_2^2} \notin Q$, then multiplicity of $\lambda_{ij}$ is $r_{ij} = 1$ and then one sensor $(D, f)$ my be sufficient to achieve $\Gamma_G$-observable of the observed systems. Now, the result bellow give information on the location of internal or boundary (pointwise and zone) $\Gamma_G$-strategic.

*4.1 Case of Zone Sensor*

We applying the previous results which are established and discussed the characterization of sensors in the case of (internal and boundary) zone sensor.

4.1.1 Internal Zone Sensor

Consider the system (16) with the output can by written by the form

$$y(t) = \int_D x(\xi_1, \xi_2, t) f(\xi_1, \xi_2) d\xi_1 d\xi_2,$$

with the zone sensor is located inside the domain $\Omega$, over the supports $D = ]\xi_1 - l_1, \xi_1 + l_1[\times]\xi_2 - l_2, \xi_2 + l_2[\subset \Omega$  (Figure 3). We have

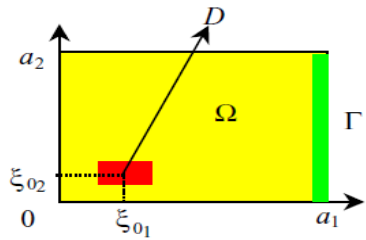

Figure 3: Location of internal zone sensor  $D$.

**Corollary 4.1:**  If  the function  $f$  is symmetric with respect to the point  $\xi_0 = (\xi_{0_1}, \xi_{0_2})$  then the sensor $(D, f)$ is  not  $\Gamma_G$-strategic if one of these conditions are satisfied:

1.    $\frac{\xi_{0_1}}{a_1} \in Q$ and $\frac{\xi_{0_2}}{a_2} \in Q$.

2.    There exist  $i_0,\ j_0 \in N$  such that $\frac{i_0 \xi_{0_1}}{a_1}$  and   $\frac{j_0 \xi_{0_2}}{a_2} \in Q$

4.1.2 Boundary Zone Sensor

Here the measurements are given by the output  $y(t) = \int_{\Gamma_0} x(\xi, t) f(\xi) d\xi$,  with  $\Gamma_0$  is an open part of  $\partial\Omega$

(Figure 4). In the case where  $\Gamma \subset \partial\Omega$  and    $f \subset L^2(\Gamma)$, the sensor  $(D, f)$  may be located on the

boundary in  $\Gamma_0 = [\eta_{0_1} - l, \eta_{0_1} + l] \times \{a_2\}$ then we have

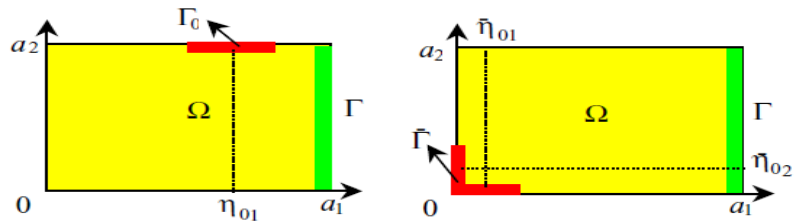

Figure 4:   Locations of boundary zone sensors $\Gamma_0, \bar{\Gamma}$.

**Corollary 4.2:**

**1. One side case:** Suppose that the sensor $(D, f)$ is located on $\Gamma_0 = [\eta_{0_1} - l, \eta_{0_1} + l] \times \{ a_2 \}$

$\subset \partial\Omega$ and $f$ is symmetric with respect to $\eta_1 = \eta_{0_1}$, then the sensor $(\Gamma_0, f)$ is not $\Gamma_G$-strategic if

$n\,\eta_{0_1} / a_1 \in Q$ for every $n, m = 1, ..., J$.

**2. Two side case:** Suppose that the sensor $(D, f)$ is located on $\bar{\Gamma} = [0, \bar{\eta}_{0_1} + l_1] \times \{0\} \bigcup \{0\}$

$\times [0, \bar{\eta}_{0_2} + l_2] \subset \partial\Omega$ and $f_{|_{\Gamma_1}}$ is symmetric with respect to $\eta_1 = \bar{\eta}_{0_1}$ and the function $f_{|_{\Gamma_2}}$ is symmetric with respect to $\eta_2 = \bar{\eta}_{0_2}$, then the sensor $(\Gamma_0, f)$ is not $\Gamma_G$-strategic if $n\,\bar{\eta}_{0_1} / a_1$ and

$m\,\bar{\eta}_{0_2} / a_2 \in Q$ for every $n, m = 1, ..., J$, where $\bar{\Gamma} = \Gamma_1 \bigcup \Gamma_2$.

This shows that the regional boundary gradient observability depends on the geometry of the sensors support and measurements function.

*4.2   Case of Pointwise Sensor*

In this subsection we discuss and characterize the sensors in the case of (internal and boundary) pointwise sensors.

### 4.2.1 Internal Pointwise Sensor

In this case the out put function is given by

$$y(t) = \int_D x(\xi_1, \xi_2, t)\delta(\xi_1 - b_1, \xi_2 - b_2)d\xi_1 d\xi_2 \qquad (19)$$

with $b = (b_1, b_2) \in \Omega$ is location of pointwise sensor    (Figure 5).

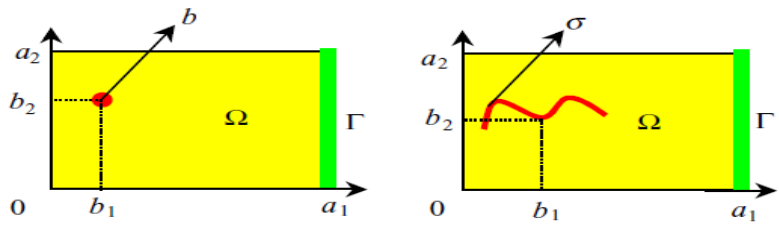

Figure 5: Locations of Internal pointwise sensors $b, \sigma$.

**Corollary 4.3:**

**1. Internal pointwise case:** If $nb_1 / a_1$ and $mb_2 / a_2 \in Q$ for every $n$, $m = 1, ..., J$, then, the sensor $(b, \delta_b)$ is not $\Gamma_G$-strategic.

**2. Filament case:** Suppose that the observation is given by the filament sensor where $\sigma = \text{Im}(\gamma)$ is symmetric with respect to the line $b = (b_1, b_2)$, if $nb_1 / a_1$ and $mb_2 / a_2 \notin Q$ for every $n$, $m = 1, ..., J$, then, the sensor $(\sigma, \delta_\sigma)$ is not $\Gamma_G$-strategic.

4.2.2 Boundary Pointwise Sensor

Here we have $b = (b_1, b_2) \in \partial\Omega$ with $b = (b_1, 0)$ or $b = (0, b_2)$ (Figure 6).

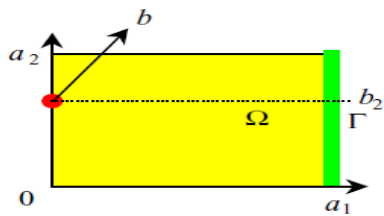

Figure 6: Location of boundary pointwise sensor $b$.

**Corollary 4.4:** The sensor $(b, \delta_b)$ is not $\Gamma_G$-strategic if $mb_2 / a_2 \in Q$ for every $m = 1, ..., J$.

This shows that there are some sensor locations to be avoided. We note that in real applications a sensor is considered as pointwise if the support area of measurement distribution is very small with respect to system domain.

**Remark 4.5:** These results can be extended to the following:

(1) Case of Neumann or mixed boundary conditions (El Jai & Pritchard 1988 ) & (Al-Saphory *et al.*2001 ).

(2) Case of disc domain with circular strategic sensor in various case of pointwise zone internal or boundary as in (Al-Saphory & El Jai 2002) & (Al-Saphory & El Jai 2001).

**5. Conclusions**

We have been characterize of regional boundary gradient strategic sensors notion in order to achieves regional boundary gradient observability. Also, we have introduced the links between the regional boundary gradient strategic sensor on Γ with a regional exactly gradient observability in ω. Thus, we have been shown that there exists a link between the exactly regional gradient observability on ω and weakly regional boundary gradient observability on Γ. Many interesting results concerning the choice of sensors structure are given and illustrated in specific situations. Various questions still open and is under consideration. For example, these result can be extended to the regional exponential general or reduced observability notions (Al-Saphory & Al-Mullah. 2015), (Al-Saphory & Jaafar  2015), observability or controllability notions for linear or non-linear (parabolic or hyperbolic) as in (Bourray *et al.* 2014), (Al-Saphory *et al.*   2010), (Ben Hadid 2012) and (Boutoulout *et al.* 2010, 2013).

**Acknowledgements.** Our thanks in advance to the editors and experts for considering this paper to publish in this esteemed journal. The authors appreciate your time and effort in reviewing the manuscript and greatly value your assistant as reviewer for the paper.

**Author's Biographies**:

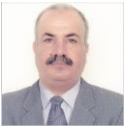

**Raheam Al-Saphory** is professor associated at the TIKRIT University in IRAQ. He received his Ph.D. degree in Control System and Analysis in (2001) from LTS/ Perpignan University France. He has a post doctorate as a researcher in 2001-2002 at LTS. Al-Saphory wrote one book and many papers in the area of control systems and analysis. Also he is a supervisor of many Ph D. and Msc. students and he was Ex-head of Department of Mathematics /College of Education for Pure Sciences, Tikrit University 2010-2011. He has visited many Centers and Scientific Departments of Bangor University/ Wales/ UK with academic staff of Iraqi Universities in 2013. He is a head of postgraduate studies committee in Mathematical Department 2014-present. Member many national and International Associations as (SMAI, SIAM, IFSA, IEEE-section IRAQ, Iraqi Mathematical and Physical Association and Iraqi Al-Khawarizemi Association).

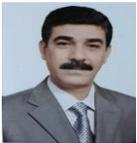

**Naseif Al-Jawari** is assistant Professor at the Al-Mustansiriyah University, IRAQ. He received his Ph.D. degree in Optimal Control Theory in (2000) from Faculty of Mathematics/ Łódź University Poland. Al-Jawari wrote many papers in the area of systems analysis and Optimization. Also he is a supervisor of many Ph D. and Msc. students. He is head of Applied Mathematics Branch/ Department of Mathematics/ College of Science/ Al-Mustansiriyah University 2012-present. Member of Iraqi Mathematical and Physical Association and Iraqi Al-Khawarizemi Association.

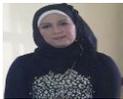

**Asmaa N  Al-Janabi** is an instructor and researcher at the Department of Mathematics/ College Science/ Al-Mustansiriyah University/ IRAQ 2012-Present. Here research area focused on Distributed Parameter Systems Analysis and Control. She is obtained here Ms.c. degree from Al-Nahrain University in (2003).